\documentclass{article}

\usepackage{amssymb}
\usepackage{amsmath}
\usepackage[margin=1.25in]{geometry}

\begin{document}

	\begin{center}
		\textbf{Large time behavior of the solution to the Cauchy problem\\
			for viscous Hamilton--Jacobi equation on infinite graphs}
		
		\medskip
		Alan A. Tedeev
	\end{center}
	
	1. Introduction
	
	We consider the large time behavior of nonnegative solutions to the Cauchy problem for a viscous Hamilton--Jacobi type equation on an infinite weighted graph of the form
	
	\begin{equation}
		\frac{\partial u}{\partial t}=\Delta _{p}u(x,t)-\frac{1}{m(x)}\sum\limits_{y\in
			V}\left\vert u(y,t)-u(x,t)\right\vert ^{q}w(x,y)\text{, }S_{T}=V\times (0,T)
		\text{,}  \tag{1.1}
	\end{equation}
	
	\begin{equation}
		u(x,0)=u_{0}(x)\geq 0\text{, }x\in V\text{.}  \tag{1.2}
	\end{equation}
	Here $V$ is the set of vertices of the graph $G(V,E,w)$ with edge set $
	E\subset V\times V$ and weight $w$.
	
	For a fixed $x_{0}$ $\in V$ let $d(x,x_{0})$ be a standard combinatorial distance in $G$ so that $d$ takes only integer values. For simplicity, we write $d(x):=d(x,x_{0})$. We assume that the graph $G$ is simple, undirected,
	infinite, connected and with locally finite degree
	
	\[
	m(x)=\sum\limits_{x\thicksim y}w(x,y),
	\]
	
	\[
	\Delta _{p}u(x,t):=\frac{1}{m(x)}\sum\limits_{x\thicksim y}\left\vert
	D_{y}u(x,t)\right\vert ^{p-2}D_{y}u(x,t)w(x,y),
	\]
	
	\[
	D_{y}u(x,t)=u(y,t)-u(x,t).
	\]
	
	In what follows we assume that
	
	\[
	q+1>p>2,\text{ }p>q.  \tag{1.3}
	\]

	The second term on the right-hand side of equation (1.1) represents a nonlinear gradient absorption. 
	
	The main goal of this paper is to determine the precise interplay between the parameters $p$, $q$, and the volume growth exponent $\alpha$ under which the total mass of the solution,
	\[
	E(t):=\sum_{x\in V}u(x,t)m(x),
	\]
	decays to zero as $t\to\infty$. We prove that there exists a critical exponent
	\[
	q^{*}=\frac{\alpha(p-1)+p}{\alpha+1},
	\]
	such that $E(t)\to 0$ whenever $q < q^*$. This condition sharply distinguishes the regime where the gradient absorption term leads to a faster decay compared to the pure $p$-Laplacian equation. We remark that the critical exponent $q^{*}$ coincides with the Fujita-type threshold obtained for related parabolic equations with gradient absorption in the continuous setting. Note that the exponent $q^*$ agrees with the critical exponent obtained in the continuous setting for degenerate parabolic equations with absorption, see \cite{andreucci2004}.
	
	To the best of our knowledge, the results presented here are new, even in the linear case $p=2$.
	
	We assume that the graph satisfies a polynomial volume growth condition: there exists $\alpha\geq 1$ and a constant $C>0$ such that for all $R\geq 1$,
	\[
	\mu(B(R)):=\sum_{x\in B(R)}m(x)\leq C R^{\alpha}. \tag{1.4}
	\]
	
	Moreover, we assume that for $q>1$ the following Hardy inequality holds:
	\begin{equation*}
		\sum_{x\in V}\frac{|u(x)|^q}{d(x)^q} m(x)
		\le C \sum_{x,y\in V} |D_y u(x)|^q w(x,y).
	\end{equation*}

	\textbf{Theorem} \textit{Let }$u(x,t)$\textit{ be a nonnegative solution of (1.1)--(1.2) in }$V\times(0,\infty)$\textit{, with }$u_{0}\in l^{1}(V)$\textit{. Suppose that the graph }$G$\textit{ satisfies the polynomial volume growth condition }$\mu(B(R))\leq C R^{\alpha}$\textit{ for some }$\alpha>0$\textit{, and let the parameters satisfy }$q+1>p>2$\textit{ and }$p>q$\textit{. Define}
	\[
	\Phi(R):=R^{(2q-p)/(q+1-p)}.
	\]
	\textit{Let }$\widetilde{R}(t)$\textit{ be the inverse function of }$\Phi(R)$\textit{, explicitly}
	\[
	\widetilde{R}(t)=t^{(q+1-p)/(2q-p)}.
	\]
	\textit{Then for all }$t>0$\textit{ the following estimate holds:}
	\[
	E(t):=\sum_{x\in V}u(x,t)m(x)\leq C\sum_{x\in V\setminus B(\widetilde{R}(t))}u_{0}(x)m(x)+C\mu(B(\widetilde{R}(t)))t^{-(p-q)/(2q-p)},
	\]
	\textit{where }$C$\textit{ is a positive constant independent of }$u_{0}$\textit{.}
	
	\textit{Consequently, using the growth condition }$\mu(B(R))\leq C R^{\alpha}$\textit{, we obtain}
	\[
	E(t)\leq C\sum_{x\in V\setminus B(\widetilde{R}(t))}u_{0}(x)m(x)+C t^{\frac{\alpha(q+1-p)-(p-q)}{2q-p}}.
	\]
	\textit{Thus }$E(t)\to 0$\textit{ as }$t\to\infty$\textit{ provided that}
	\[
	q<\frac{\alpha(p-1)+p}{\alpha+1}.
	\]
	
	\medskip
	
	The study of large time behavior of solutions to nonlinear parabolic equations on graphs has attracted considerable attention in recent years. Such problems arise naturally in the analysis of diffusion processes on discrete structures and have important applications in probability theory, mathematical physics, and network analysis.
	
	In the linear case ($p=2$), the asymptotic behavior of solutions and related properties of the discrete Laplacian have been extensively investigated; see, for instance, \cite{keller2010, grigoryan, chung, coulhon}. These works provide a detailed understanding of heat kernel estimates, spectral properties, and decay behavior on graphs with various geometric structures.
	
	The nonlinear case, involving the discrete $p$-Laplacian, is significantly more delicate due to the lack of linear structure. Problems related to existence, extinction, and long time behavior of solutions have been studied in \cite{mugnolo, hua, xin, lee, chung2014}. In particular, the interaction between diffusion and nonlinear absorption terms plays a crucial role in determining the qualitative properties of solutions.
	
	Recent works have further investigated decay properties and asymptotic estimates for nonlinear equations on infinite graphs; see \cite{andreucci2020} and the preprints \cite{tedeev_decay, tedeev_density}. However, the influence of gradient absorption terms of Hamilton--Jacobi type on the large time behavior of solutions remains less understood.
	
	The present paper contributes to this direction by establishing sharp conditions on the parameters that guarantee decay of the total mass and by identifying the critical exponent $q^*$ separating different asymptotic regimes. The proof is based on a combination of Hardy and Hölder inequalities, the use of cut-off functions adapted to the combinatorial distance, and a careful choice of the scaling parameter $R(t)=t^{(q+1-p)/(2q-p)}$.

	\newpage
	2. Proof of Theorem
	
	Denote
	
	\[
	E(t):=\sum\limits_{x\in V}u(x,t)m(x).
	\]
	
	We have
	
	\begin{equation}
		E(t)=\sum\limits_{x\in B(R)}u(x,t)m(x)+\sum\limits_{V\diagdown
			B(R)}u(x,t)m(x):=I_{1}+I_{2}.  \tag{2.1}
	\end{equation}
	
	By applying the Hölder and Hardy inequalities, we get
	
	\[
	I_{1}\leq \left( \sum\limits_{x\in B(R)}\frac{u^{q}(x,\tau )}{d(x)^{q}}
	m(x)\right) ^{1/q}\left( \sum\limits_{x\in B(R)}d(x)^{q/(q-1)}m(x)\right)
	^{(q-1)/q}
	\]
	
	\[
	\leq C\left( \sum\limits_{x,y\in V}\left\vert D_{y}u(x,\tau )\right\vert
	^{q}w(x,y)\right) ^{1/q}\left( \mu (B(R))R^{q/(q-1)}\right) ^{(q-1)/q}
	\]
	
	\begin{equation}
		=C\left( -\frac{\partial }{\partial \tau }E(\tau )\right) ^{1/q}\left( \mu
		(B(R))R^{q/(q-1)}\right) ^{(q-1)/q}.  \tag{2.2}
	\end{equation}
	
	\bigskip Next, let $\zeta (d(x))$ be a cut-off function of a ball $B_{R},$
	which we define as follows $\zeta (d(x))=1$ when $d(x)\leq 2R$, $\zeta
	(d(x))=0$ for $d(x)\geq 3R$ and for $2R\leq d(x)\leq 3R$:
	
	\[
	\zeta =\frac{d(x)-2R}{R}\text{.}
	\]
	Multiplying both sides of the equation (1.1) by $\zeta _{R}^{s}(d(x))$, $
	s>p$, and using discrete integration by parts we get
	
	\[
	\sum\limits_{x\in V}u(x,t)\zeta
	^{s}m(x)+\int\limits_{0}^{t}\sum\limits_{x,y\in V}\left\vert D_{y}u(x,\tau
	)\right\vert ^{q}\zeta ^{s}(x)w(x,y)d\tau
	\]
	
	\begin{equation}
		=-\frac{s}{2}\int\limits_{0}^{t}\sum\limits_{x,y\in V}\zeta
		^{s-1}\left\vert D_{y}u\right\vert ^{p-2}D_{y}uD_{y}\zeta
		w(x,y)+\sum\limits_{x\in V}u_{0}(x)\zeta _{R}^{s}m(x):=-\frac{s}{2}
		J_{1}+J_{2}\text{.}  \tag{2.3}
	\end{equation}
	
	\bigskip Applying the Young inequality, we obtain
	
	\[
	J_{1}:=\int\limits_{0}^{t}\sum\limits_{x,y\in V}\zeta ^{s-1}\left\vert
	D_{y}u\right\vert ^{p-2}D_{y}uD_{y}\zeta w(x,y)
	\]
	\[
	\leq \int\limits_{0}^{t}\sum\limits_{x,y\in V}\zeta ^{s-1}(x)\left\vert
	D_{y}u\right\vert ^{p-1}\left\vert D_{y}\zeta \right\vert w(x,y)
	\]
	\[
	=\int\limits_{0}^{t}\sum\limits_{x,y\in V}\zeta ^{s-1}(y)\left\vert
	D_{y}u\right\vert ^{p-1}\left\vert D_{y}\zeta \right\vert w(x,y).
	\]
	Thus,
	
	\[
	J_{1}\leq \int\limits_{0}^{t}\sum\limits_{x,y\in V}\frac{1}{2}\left( \zeta
	^{s-1}(y)+\zeta ^{s-1}(x)\right) \left\vert D_{y}u\right\vert
	^{p-1}\left\vert D_{y}\zeta \right\vert w(x,y):=\mathcal{E}(t,R)\text{.}
	\]
	
	\[
	\leq C\int\limits_{0}^{t}\sum\limits_{x,y\in V}\zeta ^{s-1}\left\vert
	D_{y}u\right\vert ^{p-1}\left\vert D_{y}\zeta \right\vert w(x,y)
	\]
	
	\[
	\leq C\varepsilon ^{q/(p-1)}\int\limits_{0}^{t}\sum\limits_{x,y\in V}\zeta
	^{s}\left\vert D_{y}u\right\vert ^{q}w(x,y)
	\]
	
	\[
	+C\varepsilon ^{-q/(q+1-p)}\int\limits_{0}^{t}\sum\limits_{x,y\in V}\zeta
	^{s-q/(q+1-p)}\left\vert D_{y}\zeta \right\vert ^{q/(q+1-p)}w(x,y)
	\]
	
	\[
	\leq C\varepsilon ^{q/(p-1)}\int\limits_{0}^{t}\sum\limits_{x,y\in V}\zeta
	^{s}\left\vert D_{y}u\right\vert ^{q}w(x,y)+Ct\mu (B(R))R^{-q/(q+1-p)}\text{.	
	}
	\]

	Therefore, from (2.3) choosing $\varepsilon $ small enough one gets
	
	\[
	\sum\limits_{x\in V}u(x,t)\zeta
	^{s}m(x)+\int\limits_{0}^{t}\sum\limits_{x,y\in V}\left\vert D_{y}u(x,\tau
	)\right\vert ^{q}\zeta ^{s}(x)w(x,y)d\tau
	\]
	
	\begin{equation}
		\leq Ct\mu (B(R))R^{-q/(q+1-p)}+C\sum\limits_{x\in V}u_{0}(x)\zeta
		_{R}^{s}m(x)\text{.}  \tag{2.4}
	\end{equation}
	
	Combining now (2.1)-(2.4), we obtain
	
	\[
	E(\tau )\leq C\left( -\frac{d}{d\tau }E(\tau )\right) ^{1/q}\left( \mu
	(B(R))R^{q/(q-1)}\right) ^{(q-1)/q}
	\]
	
	\begin{equation}
		+Ct\mu (B(R))R^{-q/(q+1-p)}+C\sum\limits_{x\in V}u_{0}(x)\zeta _{R}^{s}m(x).
		\tag{2.5}
	\end{equation}
	
	Hence we deduce
	
	\[
	F(\tau ):=F(\tau ,t,R)=E(\tau )-\left( Ct\mu
	(B(R))R^{-q/(q+1-p)}+C\sum\limits_{x\in V}u_{0}(x)\zeta _{R}^{s}m(x)\right)
	\]
	
	\begin{equation}
		\leq C\left( -\frac{d}{d\tau }F(\tau )\right) ^{1/q}\left( \mu
		(B(R))R^{q/(q-1)}\right) ^{(q-1)/q}\text{.}  \tag{2.6}
	\end{equation}
	
	Using the volume growth condition $\mu(B(R)) \le C R^\alpha$, we obtain
	\[
	\left( \mu(B(R)) R^{\frac{q}{q-1}} \right)^{\frac{q-1}{q}}
	\le C R^{\frac{\alpha(q-1)}{q} + 1}.
	\]
	Therefore, inequality (2.6) yields
	\[
	\frac{dF}{d\tau} \le - C F^q R^{-\left(\frac{\alpha(q-1)}{q} + 1\right)q}
	= - C F^q R^{-\alpha(q-1) - q}.
	\]
	Equivalently,
	\[
	\frac{dF}{F^q} \le - C R^{-\alpha(q-1) - q} \, d\tau.
	\]
	
	Hence integrating this between $0$ and $t$, we obtain
	
	\[
	F(t)\leq C\mu (B(R))\left( R^{q}t^{-1}\right) ^{1/(q-1)}\text{.}
	\]
	
	Thus,
	
	\begin{equation}
		E(t)\leq Ct\mu (B(R))\left( R^{-q/(q+1-p)}+\left( R^{q}t^{-1}\right)
		^{1/(q-1)}\right) +C\sum\limits_{x\in V}u_{0}(x)\zeta _{R}^{s}m(x).
		\tag{2.7}
	\end{equation}
	
	Finally, choosing the free parameter
	\[
	R(t)=t^{(q+1-p)/(2q-p)},
	\]
	we arrive at
	
	\[
	E(t)\leq C\sum_{x\in V\setminus B(R(t))}u_{0}(x)m(x)+C\mu(B(R(t)))t^{-(p-q)/(2q-p)}.
	\]
	
	As required. Theorem is proved. $\square $
	
	\bigskip

	The following examples illustrate the scope of the main result. \\
	\noindent\textbf{Examples.}
	
	\begin{itemize}
		
		\item Let $G=\mathbb{Z}^{N}$ with standard weights. Then $m(x)=2N$, and the volume growth exponent is $\alpha = N$, i.e., $\mu(B(R)) \sim C R^{N}$. Hence, the critical exponent becomes
		\[
		q^{*}=\frac{N(p-1)+p}{N+1}.
		\]
		For $p=2$, this gives $q^{*}=2$, which is in agreement with the scaling properties of viscous Hamilton--Jacobi equations with gradient absorption in the continuous setting.
		
		\item If the graph has exponential volume growth, the condition $\mu(B(R))\leq CR^{\alpha}$ fails for any finite $\alpha$, and the above scaling argument does not apply. In this case, different techniques are required to study the large time behavior of solutions.
		
		\item Let $G$ be a homogeneous tree (regular tree) of degree $k\geq 3$. Then
		\[
		\mu(B(R))\sim C(k-1)^{R},
		\]
		so the volume growth is exponential. Hence, the assumptions of the theorem are not satisfied, and the decay behavior may differ significantly from the polynomial growth case.
		
		\item Let $G$ be a graph with polynomial volume growth of order $D>0$, for instance a Cayley graph of a finitely generated nilpotent group. Then
		\[
		\mu(B(R))\sim CR^{D},
		\]
		and the theorem applies with $\alpha=D$. In particular, the critical exponent becomes
		\[
		q^{*}=\frac{D(p-1)+p}{D+1}.
		\]
		This shows that the result extends naturally to a broad class of graphs beyond $\mathbb{Z}^{N}$.
		
	\end{itemize}

	\bigskip
	
\end{document}